\documentclass[12pt,fleqn]{article}
\usepackage{amsmath,amssymb,amsthm,enumerate}%,backref,cite
\usepackage[cp1251]{inputenc}
\usepackage[ukrainian]{babel}
\usepackage[T2A]{fontenc}
\usepackage[all]{xy}

\begin{document}
УДК 513.83; 517.5
\bigskip

\textbf{Ю.~Б.~Зелінський, М.~В.~Стефанчук} (Yu.~B.~Zelinskii,
M.~V.~Stefanchuk), Ін-т математики НАН України, Київ.
\bigskip

\textbf{УЗАГАЛЬНЕННЯ ЗАДАЧІ ПРО ТІНЬ}

(\textbf{GENERALIZATIONS OF THE SHADOW PROBLEM})
\bigskip

\medskip
We received a solution of the shadow problem in $\mathbb{R}^{n}$ for
a family of sets, constructing from any convex domain having
nonempty interior with the help of parallel translations and
homotheties. We find a number of balls with centers on the sphere,
sufficient for giving a shadow in \emph{n}-dimensional complex
(hypercomplex) space.
\medskip

 У роботі отримано розв'язок задачі про тінь в \emph{n}-вимірному
евклідовому просторі $\mathbb{R}^{n}$ для сім'ї множин, отриманих з
довільної опуклої множини з непорожньою внутрішністю за допомогою
паралельних перенесень та гомотетій. Крім цього досліджено, яка
кількість куль з центрами на сфері достатня для створення тіні в
\emph{n}-вимірному комплексному (гіперкомплексному) просторі.
\medskip

В работе получено решение задачи о тени в \emph{n}-мерном евклидовом
пространстве $\mathbb{R}^{n}$ для семейства множеств, полученных с
произвольного выпуклого множества с непустой внутренностью с помощью
параллельных переносов и гомотетий. Кроме этого исследовано, какое
количество шаров с центрами на сфере достаточно для создания тени в
\emph{n}-мерном комплексном (гиперкомплексном) пространстве.

\newpage

\textbf{1. Вступ.} Метою роботи є розв'язок задачі про тінь для
сім'ї множин, отриманих із довільної опуклої множини з непорожньою
внутрішністю за допомогою групи перетворень~-- паралельних
перенесень та гомотетій, в \emph{n}-вимірному евклідовому просторі
$\mathbb{R}^{n}$. Це еквівалентно знаходженню умов, які забезпечують
належність точки узагальнено опуклій оболонці цієї сім'ї множин.

\textbf{Означення 1.} Множина $B\subset\mathbb{R^{\emph{n}}}$
називається \emph{m}-опуклою відносно точки
$\emph{x}\in\mathbb{R^{\emph{n}}}\setminus\emph{B}$, якщо знайдеться
\emph{m}-вимірна площина \emph{L}, така, що $\emph{x}\in\emph{L}$ і
$\emph{L}\bigcap\emph{B}=\emptyset$.

\textbf{Означення 2.} Множина $B$ називається \emph{m}-опуклою, якщо
вона \emph{m}-опукла відносно кожної точки
$\emph{x}\in\mathbb{R^{\emph{n}}}\setminus\emph{B}$.

Ці означення задовольняють аксіому опуклості~-- перетин
\emph{m}-опуклих множин буде \emph{m}-опуклою множиною. Тому для
довільної множини $B\subset\mathbb{R^{\emph{n}}}$ можна розглядати
мінімальну \emph{m}-опуклу множину, яка містить $B$, і називати її
\emph{m}-опуклою оболонкою множини $B$. [1-3]

Частковим випадком належності точки 1-оболонці об'єднання сім'ї куль
є задача про тінь, сформульована Худайбергановим [4], [5], [7].

\textbf{Задача} (про тінь). Яка мінімальна кількість попарно
неперетинних замкнених куль з центрами на сфері $\emph{S}^{n-1}$ та
радіусамии, меншими від радіуса сфери, достатня для того, щоб
довільна пряма, яка проходить через центр сфери, перетинала хоча б
одну з цих куль?

Худайберганов довів, що при $n=2$ двох куль достатньо для створення
тіні [7]. У роботах [1], [2], отримано повний розв'язок цієї задачі.
Показано, що при $n>2$ $(n+1)$-ї кулі необхідно і достатньо для
того, щоб центр сфери належав 1-опуклій оболонці сім'ї куль.

\textbf{2. 1-опуклість для центра сфери.} Нехай в \emph{n}-вимірному
евклідовому просторі $\mathbb{R}^{n}$, $n\ge2$, задана опукла
множина з непорожньою внутрішністю. Із даної множини за допомогою
групи геометричних перетворень отримана сім'я попарно неперетинних
замкнених множин. Постає проблема~-- скільки (найменша кількість)
елементів цієї сім'ї достатньо для того, щоб вибрана точка
$x\in\mathbb{R}^{n}$ належила 1-опуклій оболонці цієї сім'ї (тобто
для того, щоб довільна пряма, яка проходить через точку \emph{х},
перетинала принаймні одну з цих множин). У роботі [3] отриманий
розв'язок цієї задачі для групи геометричних перетворень, яка
складається з рухів та гомотетій опуклої множини з непорожньою
внутрішністю. Показано, що для цього необхідно і достатньо \emph{n}
елементів сім'ї. Розглянемо аналогічну задачу для сім'ї множин,
отриманих з опуклої множини з непорожньою внутрішністю за допомогою
паралельних перенесень та гомотетій.

\textbf{Теорема 1.} Для того, щоб вибрана точка в \emph{n}-вимірному
евклідовому просторі при $n\ge 2$ належала 1-опуклій оболонці сім'ї
попарно неперетинних замкнених множин, отриманих із заданої опуклої
множини з непорожньою внутрішністю за допомогою групи перетворень,
яка складається з паралельних перенесень та гомотетій, необхідно і
достатньо \emph{n} елементів цієї сім'ї.

\emph{Доведення.} Нехай \emph{В}~-- задана опукла множина з
непорожньо внутрішністю і $O\in\mathbb{R}^{n}$~-- довільна точка.
Виконаємо паралельне перенесення множини \emph{В} так, щоб точка
\emph{О} належала внутрішності множини \emph{В}. Локально межу
опуклої множини можна задати опуклою функцією [8], яка
диференційовна майже всюди на області задання [6]. Тому в будь-якому
околі довільної точки межі множини \emph{В} знайдеться точка, в якій
існує дотична гіперплощина. Опишемо навколо множини \emph{В}
многогранник з якомога меншою кількістю граней, грані якого є
дотичними гіперплощинами до множини \emph{В}. В залежності від виду
множини \emph{В} кількість граней многогранника буде лежати в межах
від $n+1$ (\emph{n}-вимірний симплекс) до $2n$ (паралелепіпед). Якщо
кількість граней більша, ніж $n+1$, то кількість пар паралельних
граней буде від 1 до \emph{n}.

Розглянемо випадок, коли навколо множини описаний паралелепіпед. З
точки \emph{О} проведемо \emph{n} променів, перпендикулярних до
кожної з граней паралелепіпеда, які не є паралельними. Нехай
$X_{1},X_{2},...,X_{n}$~-- точки перетину променів з гранями
паралелепіпеда, а $Y_{1},Y_{2},...,Y_{n}$~-- точки перетину променів
з множиною \emph{В}, інколи ці точки відповідно співпадають. На
кожному з променів $OX_{1},OX_{2},...,OX_{n}$ в досить малих околах
точок $Y_{1},Y_{2},...,Y_{n}$ виберемо точки $Z_{1},Z_{2},...,Z_{n}$
так, щоб точки $Y_{i}$ були внутрішніми точками відрізків $OZ_{i}$,
$i=1,...,n$. Виконаємо паралельні перенесення множини \emph{В} на
вектори
$\overrightarrow{Z_{1}O},\overrightarrow{Z_{2}O},...,\overrightarrow{Z_{n}O}$.
Отримаємо множини $B_{1},B_{2},...,B_{n}$. Однак утворені таким
способом множини можуть перетинатися.

Нехай $r_{1}$~-- радіус мінімального кола з центром в точці
\emph{О}, яке містить кожну з множин $B_{1},B_{2},...,B_{n}$. А
$r_{2}$~-- радіус максимального кола з центром в цій точці,
внутрішність якого не перетинається з жодною з утворених множин і
яке дотикається принаймні до однієї з множин
$B_{1},B_{2},...,B_{n}$. Виконаємо гомотетичні перетворення множин
$B_{2},...,B_{n}$ з коефіцієнтами відповідно
$$k_{1}=\frac{r_{2}}{r_{1}}, k_{2}={(\frac{r_{2}}{r_{1}})}^{2}={k_{1}}^{2},...,k_{n-1}={(\frac{r_{2}}{r_{1}})}^{n-1}={k_{1}}^{n-1}.$$
Отримаємо множини $B^{'}_{2},...,B^{'}_{n}$. Утворені таким способом
множини $B_{1},B^{'}_{2},...,B^{'}_{n}$ не перетинаються. Легко
переконатися, що двосторонній конус, заданий множинами
$B_{1},B^{'}_{2},...,B^{'}_{n}$, містить всі точки простору
$\mathbb{R}^{n}$. Таким чином, \emph{n} опуклих множин достатньо для
створення тіні. З леми 1 [1] випливає, що $(n-1)$-ї опуклої множини
замало для створення тіні. Тому для створення тіні необхідно і
достатньо \emph{n} опуклих множин. Легко переконатися, що кількість
опуклих множин залишиться такою самою, якщо навколо множини \emph{В}
описаний многогранник з меншою кількістю граней. Теорема доведена.

\textbf{3. 1-опуклість для внутрішності сфери.} Узагальнимо задачу
про тінь для куль з центрами на сфері. Скільки найменше попарно
неперетинних відкритих (замкнених) куль з центрами на сфері
$S^{n-1}$ та радіусами, меншими (які не перевищують) від радіуса
сфери, забезпечить належність внутрішності сфери 1-опуклій оболонці
сім'ї куль?

Розглянемо випадок $n=2$. Впишемо в коло правильний трикутник.
Побудуємо три відкриті (замкнені) круги з центрами у вершинах
трикутника та радіусами, які рівні половині сторони трикутника.
Бачимо, що довільна пряма, що проходить через будь-яку точку
внутрішності кола, перетне хоча б один з цих кругів. Тобто
внутрішність кола належить 1-опуклій оболонці цих кругів. Але при
цьому кулі попарно дотикаються. Для уникнення цього для досить
малого $\varepsilon$ розглянемо три круги радіусів $c+\varepsilon,
c-\varepsilon/2, c-\varepsilon/2^{2}$, де \emph{c}~-- половина
довжини сторони трикутника. Розмістимо круги так, щоб вони попарно
дотикалися, а їхні центри утворювали трикутник, який мало
відрізняється від правильного. Через центри цих кругів проходить
єдине коло, внутрішність якого належить 1-оболонці цих кругів.
Внутрішності кругів з центрами у вершинах трикутника утворюють сім'ю
з трьох відкритих куль, для якої внутрішність кола, описаного
навколо трикутника, належить 1-оболонці цієї сім'ї. Якщо круги
замкнені, то в силу неперервності, трохи зменшивши їхні радіуси,
отримаємо, що трьох замкнених куль достатньо для створення тіні для
внутрішності кола, описаного навколо трикутника. Тому справедива
наступна теорема.

\textbf{Теорема 2.} Для того, щоб внутрішність кола належала
1-опуклій оболонці сім'ї попарно неперетинних відкритих (замкнених)
кругів з центрами на колі та радіусами, меншими від радіуса кола,
необхідно і достатньо трьох кругів.

\textbf{Зауваження 1.} У випадку $n>2$ методика, застосована при
доведенні теореми 2, не проходить. При $n=3$ впишемо у двовимірну
сферу правильний тривимірний симплекс та розмістимо у його вершинах
чотири кулі з радіусами, рівними половині ребра симплекса. Тоді
через кожну з точок, яка є серединою ребра симплекса, можна провести
пряму, яка не буде перетинати жодної з цих куль.

\textbf{Зауваження 2.} Аналогічно до міркувань, які ми
використовували при доведенні теореми 2, можна показати, що при
$n>2$ існує скінченна кількість куль з центрами на сфері, для яких
довільна точка внутрішності сфери належить їх 1-опуклій оболонці.
Але знайти мінімальну кількість цих куль поки що не вдалося.

\textbf{Зауваження 3.} $(n+1)$-ї кулі достатньо, щоб усі точки, які
містяться у кулі, обмеженій сферою, належали до $(n-1)$-опуклої
оболонки системи куль. Для правильного симплекса візьмемо кулі з
центрами у його вершинах та радіусами, рівними половині висоти
симплекса. Тоді опукла оболонка цієї системи співпадає з її
$(n-1)$-опуклою оболонкою та містить сферу.

\textbf{4. Задача про тінь для 1-півопуклості.} Розглянемо більш
загальні щодо попередніх означень об'єкти.

\textbf{Означення 3.} Множина $B\subset\mathbb{R^{\emph{n}}}$
називається \emph{m}-півопуклою відносно точки
$\emph{x}\in\mathbb{R^{\emph{n}}}\setminus\emph{B}$, якщо знайдеться
\emph{m}-вимірна півплощина \emph{L}, така, що $\emph{x}\in\emph{L}$
і $\emph{L}\bigcap\emph{B}=\emptyset$.

\textbf{Означення 4.} Множина $B$ називається \emph{m}-півопуклою,
якщо вона \emph{m}-півопукла відносно кожної точки
$\emph{x}\in\mathbb{R^{\emph{n}}}\setminus\emph{B}$.

Ці означення також задовольняють аксіому опуклості. Тому для
довільної множини $B\subset\mathbb{R^{\emph{n}}}$ можна розглядати
\emph{m}-півопуклу оболонку цієї множини. [1-3]

Розглянемо аналог задачі про тінь для напівопуклості. Яка найменша
кількість  попарно неперетинних замкнених (відкритих) множин,
отриманих із заданої опуклої множини з непорожньою внутрішністю за
допомогою паралельних перенесень та гомотетій, достатня для того,
щоб вибрана точка $x\in\mathbb{R}^{n}$ належила 1-півопуклій
оболонці цієї сім'ї (тобто для того, щоб довільнй промінь, який
проходить через цю точку, перетинав принаймні одну з цих множин). У
роботі [3] отриманий розв'язок цієї задачі для групи геометричних
перетворень, яка складається з рухів та гомотетій опуклої множини з
непорожньою внутрішністю (показано, що для цього необхідно і
достатньо $(n+1)$-ї множини).

\textbf{Теорема 3.} Для того, щоб вибрана точка в \emph{n}-вимірному
евклідовому просторі при $n\ge 2$ належала 1-напівопуклій оболонці
сім'ї попарно неперетинних замкнених множин, отриманих із заданої
опуклої множини з непорожньою внутрішністю за допомогою групи
перетворень, яка складається з паралельних перенесень та гомотетій,
необхідно і достатньо 2\emph{n} елементів цієї сім'ї.

\emph{Доведення.} Як і у доведенні теореми 1 опишемо навколо опуклої
множини \emph{В} многогранник з якомога меншою кількістю граней. У
внутрішності множини \emph{В} виберемо довільну точку \emph{О}. З
цієї точки проведемо промені, перпендикулярні до кожної грані
многогранника. В залежності від виду многогранника кількість таких
променів буде лежати в межах від $n+1$ до $2n$. Нехай
$X_{1},X_{2},...,X_{n+1},...,X_{2n}$~-- точки перетину променів з
гранями многогранника, а $Y_{1},Y_{2},...,Y_{n+1},...,Y_{2n}$~--
точки перетину променів з множиною \emph{В}, інколи ці точки
відповідно співпадають. На кожному з променів
$OX_{1},OX_{2},...,OX_{n+1},...,OX_{2n}$ в досить малих околах точок
$Y_{1},Y_{2},...,Y_{n+1},...,Y_{2n}$ виберемо точки
$Z_{1},Z_{2},...,Z_{n+1},...,Z_{2n}$ так, щоб точки $Y_{i}$ були
внутрішніми точками відрізків $OZ_{i}$, $i=1,...,n+1,...,2n$ (в
залежності від виду многогранника кількість цих променів буде від
$n+1$ до $2n$). Виконаємо паралельні перенесення множини \emph{В} на
вектори
$\overrightarrow{Z_{1}O},\overrightarrow{Z_{2}O},...,\overrightarrow{Z_{n+1}O},...,\overrightarrow{Z_{2n}O}$.
Отримаємо множини $B_{1},B_{2},...,B_{n+1},...,B_{2n}$. Однак
утворені таким способом множини можуть перетинатися.

Нехай $r_{1}$~-- радіус мінімального кола з центром в точці
\emph{О}, яке містить кожну з множин
$B_{1},B_{2},...,B_{n+1},...,B_{2n}$. А $r_{2}$~-- радіус
максимального кола з центром в цій точці, внутрішність якого не
перетинається з жодною з утворених множин і яке дотикається
принаймні до однієї з множин $B_{1},B_{2},...,B_{n+1},...,B_{2n}$.
Виконаємо гомотетичні перетворення множин
$B_{2},...,B_{n+1},...,B_{2n}$ з коефіцієнтами відповідно
$$k_{1}=\frac{r_{2}}{r_{1}}, k_{2}={(\frac{r_{2}}{r_{1}})}^{2}={k_{1}}^{2},...,k_{n}={(\frac{r_{2}}{r_{1}})}^{n}=$$
$$={k_{1}}^{n},...,k_{2n-1}={(\frac{r_{2}}{r_{1}})}^{2n-1}={k_{1}}^{2n-1}.$$
Отримаємо множини $B^{'}_{2},...,B^{'}_{n+1},...,B^{'}_{2n}$.
Утворені таким способом множини $B_{1},B^{'}_{2},...,$
$B^{'}_{n+1},...,B^{'}_{2n}$ не перетинаються. Легко переконатися,
що конус, заданий множинами
$B_{1},B^{'}_{2},...,B^{'}_{n+1},...,B^{'}_{2n}$, містить всі точки
простору $\mathbb{R}^{n}$. Таким чином, якщо навколо множини
\emph{В} можна описати \emph{n}-вимірний симплекс, то $(n+1)$-ї
опуклої множини достатньо для створення тіні, якщо паралелепіпед~--
то для створення тіні достатньо $2n$ множин, в іншому випадку для
створення тіні достатньо від $n+2$ до $2n-1$ опуклих множин. З
теореми 4 [1] випливає, що $n$ опуклих множин замало для створення
тіні. Тому для створення тіні необхідно і достатньо $2n$ опуклих
множин. Теорема доведена.

\textbf{Зауваження 4.} Узагальнимо задачу про тінь для півопуклості.
Яка найменша кількість попарно неперетинних відкритих (замкнених)
куль з центрами на сфері $S^{n-1}$ та радіусами, меншими (які не
перевищують) від радіуса сфери, забезпечить належність внутрішності
сфери 1-півопуклій оболонці сім'ї куль? При $n>1$ не існує
скінченної кількості куль для створення тіні у довільній точці
внутрішності сфери. Візьмемо дві кулі, центри яких найближчі. На
відрізку, що з'єднує центри цих куль знайдеться точка, яка не
належить їх об'єднанню, та існує промінь з початком у будь-якій
точці внутрішності сфери, який проходить через дану точку, що не
перетинає систему куль.

\textbf{5. Задача про тінь в комплексному та гіперкомплексному
просторах.} Розглянемо \emph{n}-вимірні комплексний та
гіперкомплексний простори.

\textbf{Означення 5.} Множина
$B\subset\mathbb{C^{\emph{n}}}(\mathbb{H^{\emph{n}}})$ називається
\emph{m}-комплексно (\emph{m}-гіперкомплексно) опуклою відносно
точки $z\in\mathbb{C^{\emph{n}}}\setminus\emph{B}$
$(\mathbb{H^{\emph{n}}}\setminus\emph{B})$, якщо знайдеться
\emph{m}-вимірна комплексна (гіперкомплексна) площина \emph{L},
така, що $\emph{z}\in\emph{L}$ і
$\emph{L}\bigcap\emph{B}=\emptyset$.

\textbf{Означення 6.} Множина $B$ називається \emph{m}-комплексно
(\emph{m}-гіперкомплексно) опуклою, якщо вона \emph{m}-комплексно
(\emph{m}-гіперкомплексно) опукла відносно кожної точки
$z\in\mathbb{C^{\emph{n}}}\setminus\emph{B}$
$(\mathbb{H^{\emph{n}}}\setminus\emph{B})$.

Для довільної множини
$B\subset\mathbb{C^{\emph{n}}}(B\subset\mathbb{C^{\emph{n}}})$ можна
розглядати мінімальну \emph{m}-комплексно (\emph{m}-гіперкомплексно)
опуклу множину, яка містить $B$, і називати її \emph{m}-комплексною
(\emph{m}-гіперкомплексною) опуклою оболонкою множини $B$.

Зелінським [3] сформульована задача про тінь в комплексному та
гіперкомплексному просторах. Яка мінімальна кількість попарно
неперетинних замкнених куль з центрами на сфері $S^{2n-1}\subset
\mathbb{C}^{n}$ $(S^{4n-1}\subset \mathbb{H}^{n})$ та радіусами,
меншими від радіуса сфери, достатня, щоб довільна комплексна
(гіперкомплексна) пряма, яка проходить через центр сфери, перетинала
хоча б одну з цих куль (тобто для того, щоб центр сфери належав
1-комплексній чи 1-гіперкомплексній оболонці цих куль)?

Встановлено, що в комплексному (гіперкомплексному) просторі при
$n=2$ для створення тіні необхідно і достатньо 2 куль. [3]
Дослідимо, скільки таких куль достатньо для створення тіні при
$n\ge3$ в цих просторах.

\textbf{Теорема 4.} Для того, щоб вибрана точка в \emph{n}-вимірному
комплексному (гіперкомплексному) евклідовому просторі
$\mathbb{C}^{n}$ $(\mathbb{H}^{n})$, $n\ge 3$, належала
1-комплексній (1-гіперкомплексній) оболонці сім'ї попарно
неперетинних відкритих (замкнених) куль з центрами на сфері
$S^{2n-1}\subset \mathbb{C}^{n}$ $(S^{4n-1}\subset \mathbb{H}^{n})$
та з радіусами, меншими від радіуса сфери, достатньо $2n$ $(4n-2)$
куль.

\emph{Доведення.} В \emph{n}-вимірному комплексному
(гіперкомплексному) евклідовому просторі $\mathbb{C}^{n}$
$(\mathbb{H}^{n})$ розглянемо сферу $S^{2n-1}$ $(S^{4n-1})$. Через
центр цієї сфери проведемо дійсну площину розмірності $2n-1$
$(4n-3)$. Вона перетне сферу $S^{2n-1}$ $(S^{4n-1})$ по сфері
$S^{2n-2}$ $(S^{4n-4})$. Через центр сфери $S^{2n-1}$ $(S^{4n-1})$
проведемо довільну комплексну (гіперкомплексну) пряму, яка не
належить площині $\mathbb{R}^{2n-1}$ $(\mathbb{R}^{4n-3})$. Ця
комплексна (гіперкомплексна) пряма має дійсну розмірність 2 (4).
Перетин даної прямої та площини $\mathbb{R}^{2n-1}$
$(\mathbb{R}^{4n-3})$ буде містити дійсну пряму. Для того щоб центр
сфери $S^{2n-2}\subset \mathbb{R}^{2n-1}$ $(S^{4n-4}\subset
\mathbb{R}^{4n-3})$ належав 1-оболонці сім'ї попарно неперетинних
відкритих (замкнених) куль з центрами на цій сфері та з радіусами,
меншими від радіуса сфери, необхідно і достатньо $2n$ $(4n-2)$ куль.
Тому $2n$ $(4n-2)$ куль достатньо, щоб центр сфери $S^{2n-1}\subset
\mathbb{C}^{n}$ $(S^{4n-1}\subset \mathbb{H}^{n})$ належав
1-комплексній (1-гіперкомплексній) оболонці сім'ї куль з
відповідними умовами. Теорема доведена.

\vskip 2mm {\small

\newpage
Автори: \textbf{Юрій Борисович Зелінський, Марія Володимирівна
Стефанчук}
\medskip

\textbf{Інститут математики НАН України}
\medskip

Телефони: \textbf{(095)-18-62-105, (098)-76-54-988}
\medskip

E-mail: \textbf{yuzelinski@gmail.com, stefanmv43@gmail.com}
\end{document}